\documentclass[10pt]{amsart}

\usepackage[margin=1in]{geometry}
\usepackage{amsmath,amssymb,mathtools}
\usepackage{chngcntr}
\usepackage{microtype}
\usepackage{xcolor}
\usepackage[colorlinks=true,linkcolor=red,citecolor=green!50!black,urlcolor=blue]{hyperref}

\counterwithout{equation}{section}

\newtheorem{theorem}{Theorem}[section]
\newtheorem{corollary}[theorem]{Corollary}
\newtheorem{lemma}[theorem]{Lemma}
\newtheorem{proposition}[theorem]{Proposition}
\theoremstyle{remark}
\newtheorem{remark}[theorem]{Remark}
\newtheorem{setting}[theorem]{Setting}

\newcommand{\N}{N}
\newcommand{\V}{V}
\newcommand{\E}{E}
\newcommand{\PG}{P_G}
\newcommand{\PH}{P_H}
\newcommand{\rhoTwo}{\rho^{\{2\}}}

\title[A $2/3$ Bound for Vizing's Conjecture]
{A $2/3$ Bound for Vizing's Conjecture}

\author{Mohsen Aliabadi}
\address{Department of Mathematics, Clayton State University,
2000 Clayton State Boulevard, Morrow, GA 30260, USA}
\email{MohsenAliabadi@clayton.edu}

\author{Elliot Krop}
\address{Department of Mathematics, Clayton State University,
2000 Clayton State Boulevard, Morrow, GA 30260, USA}
\email{elliotkrop@clayton.edu}
\date{August 26, 2026}

\begin{document}

\begin{abstract}
Vizing's conjecture, dating back to 1963, asserts that
\[
\gamma(G\mathbin{\square}H)
\geq
\gamma(G)\gamma(H)
\]
for all finite graphs $G$ and $H$, where $\gamma$ denotes the domination
number and $\square$ denotes the Cartesian product. In 2000, Clark and Suen
proved the universal bound
\[
\gamma(G\mathbin{\square}H)
\geq
\frac{1}{2}\gamma(G)\gamma(H).
\]
Recently, Steiner obtained the first constant-factor improvement of the
Clark--Suen bound, proving that
\[
\gamma(G\mathbin{\square}H)
\geq
\frac{5+\sqrt{73}}{24}\gamma(G)\gamma(H)
\approx
0.5643\,\gamma(G)\gamma(H).
\]
In this paper, we further improve the universal constant by proving that
\[
\gamma(G\mathbin{\square}H)
\geq
\frac{2}{3}\gamma(G)\gamma(H)
\]
for all finite graphs $G$ and $H$. Thus, we raise the best known universal
constant in the approximate form of Vizing's conjecture from
$(5+\sqrt{73})/24$ to $2/3$.
\end{abstract}

\maketitle

\maketitle

\section{Introduction}

All graphs in this paper are finite, simple, and have at least one vertex.
For a graph $G$ and a vertex $v$, $\N_G(v)$ denotes the open neighborhood of
$v$ and $\N_G[v]=\N_G(v)\cup\{v\}$ its closed neighborhood. A set
$D\subseteq \V(G)$ dominates a vertex $v$ if $\N_G[v]\cap D\neq\varnothing$,
and dominates a set $S\subseteq \V(G)$ if it dominates every vertex of $S$.
The domination number $\gamma(G)$ is the minimum size of a set dominating
$\V(G)$. The Cartesian product $G\mathbin{\square}H$ has vertex set
$\V(G)\times\V(H)$; two distinct vertices $(g,h)$ and $(g',h')$ are adjacent
if either $g=g'$ and $hh'\in\E(H)$, or $h=h'$ and $gg'\in\E(G)$.

\medskip
\noindent\textbf{Vizing's conjecture \cite{Vizing1963,Vizing1968}.}
For all graphs $G$ and $H$,
\[
 \gamma(G\mathbin{\square}H)\geq \gamma(G)\gamma(H).
\]

The conjecture is one of the central open problems in domination theory; see
the survey \cite{BresarEtAl2012} and the monographs \cite{Haynes2023,Haynes1998}
for its history. For restricted graph classes, stronger bounds are known. In particular,
the second author of the present paper ~\cite{Krop2017} proved the $2/3$-bound when one factor is claw-free;
this was subsequently improved to $3/4$ for claw-free graphs by
Bre\v{s}ar and Henning~\cite{BresarHenning2020}. Our result establishes
the $2/3$-bound for arbitrary graphs.
Among results valid for all graphs, the classical bound of
Clark and Suen \cite{ClarkSuen2000},
\begin{equation}
 \gamma(G\mathbin{\square}H)\geq \frac12\gamma(G)\gamma(H),
 \label{eq:clark-suen}
\end{equation}
was improved in lower-order terms by Suen and Tarr \cite{SuenTarr2012} and by
Zerbib \cite{Zerbib2019}. Bre\v{s}ar \cite{Bresar2017} proved
\begin{equation}
 \gamma(G\mathbin{\square}H)
 \geq \frac{2\gamma(G)-\rho(G)}{3}\gamma(H),
 \label{eq:bresar}
\end{equation}
where $\rho(G)$ is the $2$-packing number of $G$ (see Section~\ref{sec:prelim}).
Very recently, Steiner \cite{Steiner2026} obtained the first constant-factor
improvement of \eqref{eq:clark-suen} valid without any additional hypothesis:
\begin{equation}
 \gamma(G\mathbin{\square}H)
 \geq \frac{5+\sqrt{73}}{24}\gamma(G)\gamma(H)
 \approx 0.5643\gamma(G)\gamma(H).
 \label{eq:steiner-main}
\end{equation}
A key ingredient of \cite{Steiner2026} is the inequality
\begin{equation}
 \gamma(G\mathbin{\square}H)
 \geq \frac{3\gamma(G)-\rhoTwo(G)}{4}\gamma(H),
 \label{eq:steiner-packing}
\end{equation}
where $\rhoTwo(G)$ is the maximum total weight of a $2$-packing function on
$G$ (Section~\ref{sec:prelim}). Steiner combines \eqref{eq:bresar} and
\eqref{eq:steiner-packing} with a bound derived from a theorem of Hou and Lu
\cite{HouLu2009}.

\medskip
\noindent\textbf{Main results.}
Inequalities \eqref{eq:bresar} and \eqref{eq:steiner-packing} are proved by
partitioning $\V(H)$ into parts contained in the closed neighborhoods of the
vertices of a minimum dominating set of $H$, and by charging each part a
penalty, respectively $\rho(G)/2$ and $\rhoTwo(G)/3$, coming from an auxiliary
estimate on the domination number of a subset of $\V(G)$. Our starting
observation is that a part equal to the full closed neighborhood $\N_H[h]$ of
a vertex $h$ can never be vertically dominated in the sense of that argument,
so such a part need not be charged. Choosing these distinguished parts to be
the closed neighborhoods of a $2$-packing of $H$, and keeping the inequalities
for the different types of parts separate rather than adding them, gives the
following theorem.

\begin{theorem}\label{thm:main-packing}
Let $G$ and $H$ be graphs, let $P\subseteq\V(H)$ be a $2$-packing of $H$, and
put
\[
 r:=|P|,\qquad W:=\V(H)\setminus\N_H[P],\qquad s:=\gamma_H(W).
\]
Then
\begin{equation}
 \gamma(G\mathbin{\square}H)
 \geq \frac{r\gamma(G)}{3}
 +\frac{(r+s)(2\gamma(G)-\rho(G))}{3}.
 \label{eq:main-packing}
\end{equation}
\end{theorem}

Here $\gamma_H(W)$ is the minimum number of vertices of $H$ dominating $W$;
in particular, $\gamma_H(\varnothing)=0$. Taking $P$ to be a maximum
$2$-packing and using $r+s\geq\gamma(H)$ gives the following consequences.

\begin{corollary}\label{cor:asymmetric}
For all graphs $G$ and $H$,
\begin{equation}
 \gamma(G\mathbin{\square}H)
 \geq
 \frac{2\gamma(G)\gamma(H)+\rho(H)\gamma(G)-\rho(G)\gamma(H)}{3},
 \label{eq:asymmetric}
\end{equation}
and hence, by symmetry,
\begin{equation}
 \gamma(G\mathbin{\square}H)
 \geq
 \frac{2\gamma(G)\gamma(H)
 +|\rho(H)\gamma(G)-\rho(G)\gamma(H)|}{3}.
 \label{eq:symmetric}
\end{equation}
\end{corollary}

\begin{theorem}\label{thm:two-thirds}
For all graphs $G$ and $H$,
\[
 \gamma(G\mathbin{\square}H)\geq \frac23\gamma(G)\gamma(H).
\]
\end{theorem}

Inequality \eqref{eq:asymmetric} exceeds Bre\v{s}ar's bound
\eqref{eq:bresar} by exactly $\rho(H)\gamma(G)/3$. Theorem~\ref{thm:two-thirds}
does not use $2$-packing functions or the bound of Hou and Lu; the only
auxiliary estimate needed is an elementary lemma implicit in
\cite{Bresar2017} (Lemma~\ref{lem:bresar-steiner} below). The same refinement
applied to Steiner's argument gives the following.

\begin{theorem}\label{thm:packing-function}
In the setting of Theorem~\ref{thm:main-packing},
\begin{equation}
 \gamma(G\mathbin{\square}H)
 \geq r\gamma(G)-\frac{r\rhoTwo(G)}8
 +\frac{s(3\gamma(G)-\rhoTwo(G))}{4}.
 \label{eq:packing-function}
\end{equation}
Consequently, for all graphs $G$ and $H$,
\begin{equation}
 \gamma(G\mathbin{\square}H)
 \geq
 \frac{2(3\gamma(G)-\rhoTwo(G))\gamma(H)
 +\rho(H)(2\gamma(G)+\rhoTwo(G))}{8}.
 \label{eq:packing-function-global}
\end{equation}
\end{theorem}

Inequality \eqref{eq:packing-function-global} exceeds Steiner's bound
\eqref{eq:steiner-packing} by
$\rho(H)(2\gamma(G)+\rhoTwo(G))/8$.

Finally, we show that the constant $2/3$ cannot be improved by combining the
inequalities of this paper with those of \cite{Steiner2026}, if the relevant
graph parameters are treated as independent quantities
(Proposition~\ref{prop:limits}). We also describe the extremal solution of the
underlying linear program, indicating where a further improvement would have
to come from.

\medskip
\noindent\textbf{Organization.}
In Section~\ref{sec:prelim} we collect definitions and auxiliary lemmas.
Section~\ref{sec:partition} sets up the partition and proves the four counting
inequalities on which everything rests. Section~\ref{sec:proofs} proves
Theorems~\ref{thm:main-packing} and \ref{thm:packing-function}, and
Section~\ref{sec:constant} derives the corollaries and
Theorem~\ref{thm:two-thirds}. Section~\ref{sec:limits} is devoted to the limits of
the method.

\section{Preliminaries}\label{sec:prelim}

For $S\subseteq\V(G)$, let $\gamma_G(S)$ be the minimum size of a set
$D\subseteq\V(G)$ dominating $S$. Thus $\gamma_G(\V(G))=\gamma(G)$ and
$\gamma_G(\varnothing)=0$. Domination is subadditive:
$\gamma_G(S\cup T)\leq\gamma_G(S)+\gamma_G(T)$.

A set $P\subseteq\V(G)$ is a $2$-packing if the closed neighborhoods of its
vertices are pairwise disjoint; equivalently,
$|\N_G[v]\cap P|\leq 1$ for every $v\in\V(G)$. The $2$-packing number
$\rho(G)$ is the maximum size of a $2$-packing. For $S\subseteq\V(G)$, let
$\rho_G(S)$ be the maximum size of a $2$-packing contained in $S$; clearly
$\rho_G(S)\leq\rho(G)$.

Following \cite{GoddardHenning2020,Steiner2026}, for an integer $k\geq0$ a
function $f:\V(G)\to\mathbb Z_{\geq0}$ is a $k$-packing function if
$f(\N_G[v])\leq k$ for every $v\in\V(G)$, where
$f(S):=\sum_{x\in S}f(x)$. The maximum total weight $f(\V(G))$ of a
$k$-packing function is denoted $\rho^{\{k\}}(G)$. The indicator function of
a $2$-packing is a $1$-packing function, so $\rho^{\{1\}}(G)=\rho(G)$.
Doubling a $1$-packing function gives a $2$-packing function, so
$\rhoTwo(G)\geq 2\rho(G)$.

\begin{lemma}\label{lem:packing-upper}
For every graph $G$,
\[
 \rho(G)\leq\gamma(G)
 \qquad\text{and}\qquad
 \rhoTwo(G)\leq2\gamma(G).
\]
\end{lemma}

\begin{proof}
Let $D$ be a minimum dominating set of $G$. If $P$ is a $2$-packing, each
$d\in D$ dominates at most one vertex of $P$, whereas $D$ dominates all of
$P$; hence $|P|\leq|D|$. If $f$ is a $2$-packing function, every vertex lies
in $\N_G[d]$ for some $d\in D$, so, as $f\geq0$,
\[
 f(\V(G))\leq\sum_{d\in D}f(\N_G[d])\leq2|D|.
\]
\end{proof}

The following estimate is Lemma 2.2 of \cite{Steiner2026}; in the weaker form
with $\rho(G)$ in place of $\rho_G(S)$ it is implicit in
\cite{Bresar2017}. We include the short proof.

\begin{lemma}[Bre\v{s}ar; Steiner]\label{lem:bresar-steiner}
For every graph $G$ and every $S\subseteq\V(G)$,
\[
 \gamma_G(S)\leq\frac{|S|+\rho_G(S)}2
 \leq\frac{|S|+\rho(G)}2.
\]
\end{lemma}

\begin{proof}
Let $B$ be the graph on vertex set $S$ in which distinct $s_1,s_2\in S$ are
adjacent if and only if $\N_G[s_1]\cap\N_G[s_2]\neq\varnothing$. Let
$\{s_1s'_1,\ldots,s_ts'_t\}$ be a maximum matching of $B$, and choose
$d_i\in\N_G[s_i]\cap\N_G[s'_i]$ for each $i$. The set
\[
 \bigl(S\setminus\{s_1,s'_1,\ldots,s_t,s'_t\}\bigr)
 \cup\{d_1,\ldots,d_t\}
\]
dominates $S$, so $\gamma_G(S)\leq|S|-t$. The unmatched vertices form an
independent set in $B$ by maximality of the matching, and hence a $2$-packing
of $G$ contained in $S$. Thus $|S|-2t\leq\rho_G(S)$. Combining these
inequalities gives
$\gamma_G(S)\leq|S|-t\leq(|S|+\rho_G(S))/2$.
\end{proof}

\begin{lemma}[Steiner \cite{Steiner2026}]\label{lem:steiner-subset}
For every graph $G$ and every $S\subseteq\V(G)$,
\[
 \gamma_G(S)\leq\frac{|S|+\rhoTwo(G)}3.
\]
\end{lemma}

Lemma~\ref{lem:steiner-subset} is used only in
Theorem~\ref{thm:packing-function}; Theorems~\ref{thm:main-packing} and
\ref{thm:two-thirds} depend on Lemma~\ref{lem:bresar-steiner} alone.

\medskip
\noindent\textbf{Products, layers, and projections.}
For $(g,h)\in\V(G\mathbin{\square}H)$,
\begin{equation}
 \N_{G\mathbin{\square}H}[(g,h)]
 =\bigl(\N_G[g]\times\{h\}\bigr)
 \cup\bigl(\{g\}\times\N_H[h]\bigr).
 \label{eq:closed-product}
\end{equation}
For $h\in\V(H)$, the set $\V(G)\times\{h\}$ is the $G$-layer at $h$; for
$g\in\V(G)$, the set $\{g\}\times\V(H)$ is the $H$-fiber at $g$. Write
$\PG$ and $\PH$ for the two coordinate projections. Notice that $\PH$ is
injective on every $H$-fiber. Throughout, $D$ denotes a fixed minimum
dominating set of $G\mathbin{\square}H$, so
$|D|=\gamma(G\mathbin{\square}H)$, and
$D_g:=D\cap(\{g\}\times\V(H))$.

\begin{lemma}\label{lem:layer-dominates}
For every $h\in\V(H)$, the set
\[
 \PG\bigl(D\cap(\V(G)\times\N_H[h])\bigr)
\]
dominates $G$. In particular,
$|D\cap(\V(G)\times\N_H[h])|\geq\gamma(G)$.
\end{lemma}

\begin{proof}
Let $g\in\V(G)$. By \eqref{eq:closed-product}, the vertex $(g,h)$ is
dominated by some $(g',h')\in D$ with either $h'=h$ and
$g'\in\N_G[g]$, or $g'=g$ and $h'\in\N_H(h)$. In both cases
$h'\in\N_H[h]$ and $g'\in\N_G[g]$.
\end{proof}

\section{The partition framework}\label{sec:partition}

\begin{setting}\label{setting:partition}
Let $G$ and $H$ be graphs and let
$P=\{h_1,\ldots,h_r\}$ be a $2$-packing of $H$, possibly empty. Put
\[
 W:=\V(H)\setminus\N_H[P],\qquad s:=\gamma_H(W),
\]
and let $Q=\{h_{r+1},\ldots,h_{r+s}\}$ be a minimum set of vertices of $H$
dominating $W$. Define a partition
$\V(H)=\pi_1\mathbin{\dot\cup}\cdots\mathbin{\dot\cup}\pi_{r+s}$ as follows.
\begin{enumerate}
\item For $i\leq r$, let $\pi_i:=\N_H[h_i]$. These sets are pairwise
disjoint because $P$ is a $2$-packing, and their union is $\N_H[P]$.
\item For $i>r$, choose $\pi_i\subseteq\N_H[h_i]\cap W$ so that
$\pi_{r+1},\ldots,\pi_{r+s}$ partition $W$. Such sets exist because $Q$
dominates $W$: assign each $w\in W$ to one index $i>r$ for which
$w\in\N_H[h_i]$.
\end{enumerate}
The parts $\pi_1,\ldots,\pi_r$ are called \emph{packing parts}, and the
remaining parts are called \emph{leftover parts}. If $s=0$, then $W$ is
empty and there are no leftover parts; if $r=0$, then $W=\V(H)$ and $Q$ is a
minimum dominating set of $H$.

For $i\in[r+s]$, let
$D_i:=D\cap(\V(G)\times\pi_i)$, and for $g\in\V(G)$ let
$C_g^i:=\{g\}\times\pi_i$. The cell $C_g^i$ is \emph{vertically dominated}
if every vertex of $C_g^i$ has a neighbor in
\[
 D\cap\bigl(\{g\}\times(\V(H)\setminus\pi_i)\bigr),
\]
and put
\[
 L_i:=\{g\in\V(G):C_g^i\text{ is vertically dominated}\}.
\]
Finally, split $D$ according to layers:
\[
 a:=|D\cap(\V(G)\times P)|,
\quad
 b:=|D\cap(\V(G)\times(\N_H[P]\setminus P))|,
\quad
 w:=|D\cap(\V(G)\times W)|.
\]
Thus
\[
 |D|=a+b+w,\qquad
 \sum_{i\leq r}|D_i|=a+b,\qquad
 \sum_{i>r}|D_i|=w.
\]
\end{setting}

The first two lemmas are the standard part of the argument of
\cite{Bresar2017,Steiner2026}, formulated for an arbitrary part.

\begin{lemma}\label{lem:horizontal}
For every $i\in[r+s]$, the set $\PG(D_i)$ dominates $\V(G)\setminus L_i$.
Consequently,
\[
 \gamma_G(\V(G)\setminus L_i)\leq|D_i|.
\]
\end{lemma}

\begin{proof}
Let $g\in\V(G)\setminus L_i$. Since $C_g^i$ is not vertically dominated,
some $(g,h)\in C_g^i$ has no neighbor in
$D\cap(\{g\}\times(\V(H)\setminus\pi_i))$. The vertex $(g,h)$ is dominated
by some $(g',h')\in D$. By \eqref{eq:closed-product}, either $g'=g$ and
$h'\in\N_H[h]$, in which case $h'\in\pi_i$ by the choice of $h$, or
$h'=h\in\pi_i$ and $g'\in\N_G[g]$. In both cases
$(g',h')\in D_i$ and $g\in\N_G[g']$.
\end{proof}

\begin{lemma}\label{lem:packing-L-empty}
For every $i\leq r$, we have $L_i=\varnothing$, and hence
$\gamma(G)\leq|D_i|$.
\end{lemma}

\begin{proof}
Let $i\leq r$ and $g\in\V(G)$. By \eqref{eq:closed-product}, the neighbors
of $(g,h_i)$ in the $H$-fiber at $g$ are the vertices $(g,h')$ with
$h'\in\N_H(h_i)\subseteq\N_H[h_i]=\pi_i$. Thus $(g,h_i)\in C_g^i$ has no
neighbor in $\{g\}\times(\V(H)\setminus\pi_i)$, so $C_g^i$ is not vertically
dominated. The second assertion follows from Lemma~\ref{lem:horizontal}.
\end{proof}

The next lemma controls the total size of the sets $L_i$ belonging to leftover
parts. Compared with \cite{Bresar2017,Steiner2026}, minimality is required only
relative to $W$, and vertices of $D$ lying above $P$ are not charged.

\begin{lemma}\label{lem:sum-L}
\[
 \sum_{i=r+1}^{r+s}|L_i|\leq |D|-a=b+w.
\]
\end{lemma}

\begin{proof}
For $g\in\V(G)$, let $m_g$ be the number of indices $i>r$ with $g\in L_i$.
Then $\sum_{i>r}|L_i|=\sum_gm_g$, and it is enough to show
\begin{equation}
 m_g\leq |D_g|-|D_g\cap(\{g\}\times P)|
 \qquad\text{for every }g\in\V(G).
 \label{eq:mg}
\end{equation}
Fix $g$ and put $X:=\PH(D_g)\setminus P$. We claim that
\[
 Y:=X\cup\{h_i:i>r,\ g\notin L_i\}
\]
dominates $W$ in $H$. Let $u\in W$, and let $i>r$ be the index with
$u\in\pi_i$. If $g\notin L_i$, then $h_i\in Y$ dominates $u$ because
$\pi_i\subseteq\N_H[h_i]$. If $g\in L_i$, then the cell $C_g^i$ is
vertically dominated, so $(g,u)$ has a neighbor $(g,h')\in D$ with
$h'\notin\pi_i$. By \eqref{eq:closed-product}, $h'\in\N_H(u)$. Since
$u\in W=\V(H)\setminus\N_H[P]$, no neighbor of $u$ lies in $P$. Hence
$h'\notin P$, so $h'\in X\subseteq Y$ and $h'$ dominates $u$.

Since $s=\gamma_H(W)$,
\[
 s\leq |Y|\leq |X|+(s-m_g),
\]
and therefore $m_g\leq|X|$. Finally, $\PH$ is injective on $D_g$, so
\[
 |X|=|D_g|-|D_g\cap(\{g\}\times P)|,
\]
which proves \eqref{eq:mg}. Summing over $g$ gives the result.
\end{proof}

We now derive the four counting inequalities.

\begin{lemma}[Packing blocks]\label{lem:packing-blocks}
\[
 a+b\geq r\gamma(G).
\]
\end{lemma}

\begin{proof}
By Lemma~\ref{lem:layer-dominates},
$|D_i|=|D\cap(\V(G)\times\N_H[h_i])|\geq\gamma(G)$ for each $i\leq r$.
Sum over $i\leq r$.
\end{proof}

\begin{lemma}[Second count on the packing blocks]\label{lem:second-count}
\begin{align}
 a+\frac b2&\geq r\left(\gamma(G)-\frac{\rho(G)}2\right),
 \label{eq:second-rho}\\
 a+\frac b3&\geq r\left(\gamma(G)-\frac{\rhoTwo(G)}3\right).
 \label{eq:second-rho2}
\end{align}
\end{lemma}

\begin{proof}
Fix $i\leq r$, write $p=h_i$, and let
\[
 A_p:=\{g\in\V(G):D\cap(\{g\}\times\N_H(p))\neq\varnothing\}.
\]
If $g\notin A_p$, then by \eqref{eq:closed-product} the vertex $(g,p)$ is
dominated by some $(g',p)\in D$ with $g'\in\N_G[g]$. Hence
$\PG(D\cap(\V(G)\times\{p\}))$ dominates $\V(G)\setminus A_p$, and
subadditivity gives
\[
 \gamma(G)
 \leq \gamma_G(\V(G)\setminus A_p)+\gamma_G(A_p)
 \leq |D\cap(\V(G)\times\{p\})|+\gamma_G(A_p).
\]
Moreover,
$|A_p|\leq|D\cap(\V(G)\times\N_H(p))|$, since distinct elements of $A_p$
are witnessed by distinct vertices of $D$. Lemma~\ref{lem:bresar-steiner}
therefore gives
\[
 \gamma(G)
 \leq |D\cap(\V(G)\times\{p\})|
 +\frac{|D\cap(\V(G)\times\N_H(p))|+\rho(G)}2.
\]
Using Lemma~\ref{lem:steiner-subset} instead gives the analogous inequality
with denominator $3$ and $\rhoTwo(G)$ in the numerator.

Now sum over $i\leq r$. The singletons $\{h_i\}$ are disjoint with union
$P$, while the open neighborhoods $\N_H(h_i)$ are pairwise disjoint, disjoint
from $P$, and have union $\N_H[P]\setminus P$. The two summed inequalities
are exactly \eqref{eq:second-rho} and \eqref{eq:second-rho2}.
\end{proof}

\begin{lemma}[Leftover parts]\label{lem:leftover}
\begin{align}
 s\gamma(G)&\leq w+\frac12\sum_{i>r}|L_i|+\frac{s\rho(G)}2,
 \label{eq:leftover-rho}\\
 s\gamma(G)&\leq w+\frac13\sum_{i>r}|L_i|+\frac{s\rhoTwo(G)}3.
 \label{eq:leftover-rho2}
\end{align}
\end{lemma}

\begin{proof}
For $i>r$, subadditivity, Lemma~\ref{lem:horizontal}, and
Lemma~\ref{lem:bresar-steiner} applied to $S=L_i$ give
\[
 \gamma(G)
 \leq\gamma_G(\V(G)\setminus L_i)+\gamma_G(L_i)
 \leq|D_i|+\frac{|L_i|+\rho(G)}2.
\]
Using Lemma~\ref{lem:steiner-subset} gives instead
\[
 \gamma(G)\leq|D_i|+\frac{|L_i|+\rhoTwo(G)}3.
\]
Sum over the $s$ indices $i>r$ and use $\sum_{i>r}|D_i|=w$.
\end{proof}

\begin{remark}\label{rem:recover-old}
Adding Lemma~\ref{lem:packing-L-empty}, summed over $i\leq r$, to
\eqref{eq:leftover-rho2}, and using the weaker estimate
$\sum_{i>r}|L_i|\leq|D|$, recovers for $r=0$ exactly the proof of
\eqref{eq:steiner-packing} in \cite{Steiner2026}; similarly,
\eqref{eq:leftover-rho} recovers \eqref{eq:bresar}. The improvements below
come from four features: the packing parts contribute $\gamma(G)$ each without
penalty; Lemma~\ref{lem:sum-L} does not charge vertices above $P$;
Lemma~\ref{lem:second-count} supplies a second count; and the inequalities of
Lemmas~\ref{lem:packing-blocks} and \ref{lem:leftover} are kept separate until
the final combination.
\end{remark}

\section{Proofs of Theorems~\ref{thm:main-packing} and
\ref{thm:packing-function}}\label{sec:proofs}

\begin{proof}[Proof of Theorem~\ref{thm:main-packing}]
Adopt Setting~\ref{setting:partition} and write
$\gamma=\gamma(G)$ and $\rho=\rho(G)$. By
\eqref{eq:leftover-rho} and Lemma~\ref{lem:sum-L},
\[
 s\gamma\leq w+\frac{b+w}{2}+\frac{s\rho}{2},
\]
and hence
\[
 w\geq\frac23s\left(\gamma-\frac\rho2\right)-\frac b3.
\]
Therefore
\begin{equation}
 |D|=a+b+w
 \geq a+\frac23b+\frac23s\left(\gamma-\frac\rho2\right).
 \label{eq:first-combination}
\end{equation}
By Lemma~\ref{lem:packing-blocks} and \eqref{eq:second-rho},
\begin{align}
 a+\frac23b
 &=\frac13(a+b)+\frac23\left(a+\frac b2\right)\notag\\
 &\geq\frac13r\gamma
 +\frac23r\left(\gamma-\frac\rho2\right).
 \label{eq:correct-combination}
\end{align}
Combining \eqref{eq:first-combination} and
\eqref{eq:correct-combination} yields
\[
 \gamma(G\mathbin{\square}H)=|D|
 \geq\frac{r\gamma}{3}
 +\frac23(r+s)\left(\gamma-\frac\rho2\right)
 =\frac{r\gamma}{3}+\frac{(r+s)(2\gamma-\rho)}3.
\]
\end{proof}

\begin{proof}[Proof of Theorem~\ref{thm:packing-function}]
Write $\gamma=\gamma(G)$ and $\rhoTwo=\rhoTwo(G)$. By
\eqref{eq:leftover-rho2} and Lemma~\ref{lem:sum-L},
\[
 s\gamma\leq w+\frac{b+w}{3}+\frac{s\rhoTwo}{3},
\]
so
\[
 w\geq\frac34s\left(\gamma-\frac{\rhoTwo}{3}\right)-\frac b4.
\]
Thus
\[
 |D|\geq a+\frac34b
 +\frac34s\left(\gamma-\frac{\rhoTwo}{3}\right).
\]
By Lemma~\ref{lem:packing-blocks} and \eqref{eq:second-rho2},
\begin{align*}
 a+\frac34b
 &=\frac58(a+b)+\frac38\left(a+\frac b3\right)\\
 &\geq\frac58r\gamma
 +\frac38r\left(\gamma-\frac{\rhoTwo}{3}\right)
 =r\gamma-\frac{r\rhoTwo}{8}.
\end{align*}
This proves \eqref{eq:packing-function}.

For \eqref{eq:packing-function-global}, let $P$ be a maximum $2$-packing of
$H$, so $r=\rho(H)$. Since $P$ dominates $\N_H[P]$ and $Q$ dominates $W$,
the set $P\cup Q$ dominates $H$. Hence
$r+s\geq\gamma(H)$ and $s\geq\gamma(H)-\rho(H)$. By
Lemma~\ref{lem:packing-upper},
$3\gamma-\rhoTwo\geq\gamma>0$, so substituting the lower bound on $s$ into
\eqref{eq:packing-function} gives
\begin{align*}
 \gamma(G\mathbin{\square}H)
 &\geq \rho(H)\gamma-\frac{\rho(H)\rhoTwo}{8}
 +\frac{(\gamma(H)-\rho(H))(3\gamma-\rhoTwo)}4\\
 &=\frac{2(3\gamma-\rhoTwo)\gamma(H)
 +\rho(H)(2\gamma+\rhoTwo)}8.
\end{align*}
\end{proof}

\begin{remark}\label{rem:choice-P}
For $P=\varnothing$, so $r=0$ and $s=\gamma(H)$,
Theorem~\ref{thm:main-packing} reduces to Bre\v{s}ar's inequality
\eqref{eq:bresar}, and Theorem~\ref{thm:packing-function} reduces to Steiner's
inequality \eqref{eq:steiner-packing}. For any $2$-packing $P$, the right-hand
side of \eqref{eq:main-packing} is at least
\[
 \frac{r\gamma(G)}3
 +\frac{\gamma(H)(2\gamma(G)-\rho(G))}3,
\]
which is increasing in $r$. Thus a maximum $2$-packing is the best choice when
nothing beyond $r+s\geq\gamma(H)$ is known about $s$, although a nonmaximum
$P$ with large $s$ can give a better bound.
\end{remark}

\section{The constant \texorpdfstring{$2/3$}{2/3}}\label{sec:constant}

\begin{proof}[Proof of Corollary~\ref{cor:asymmetric}]
Apply Theorem~\ref{thm:main-packing} to a maximum $2$-packing $P$ of $H$.
Then $r=\rho(H)$ and, as observed above, $r+s\geq\gamma(H)$. By
Lemma~\ref{lem:packing-upper},
$2\gamma(G)-\rho(G)\geq\gamma(G)>0$, so \eqref{eq:main-packing} gives
\begin{align*}
 \gamma(G\mathbin{\square}H)
 &\geq\frac{\rho(H)\gamma(G)}3
 +\frac{\gamma(H)(2\gamma(G)-\rho(G))}3\\
 &=\frac{2\gamma(G)\gamma(H)+\rho(H)\gamma(G)-\rho(G)\gamma(H)}3,
\end{align*}
which is \eqref{eq:asymmetric}. Since
$G\mathbin{\square}H\cong H\mathbin{\square}G$, the same inequality holds
with the roles of $G$ and $H$ interchanged. Taking the larger of the two bounds
gives \eqref{eq:symmetric}.
\end{proof}

\begin{proof}[Proof of Theorem~\ref{thm:two-thirds}]
This follows immediately from \eqref{eq:symmetric}, since the absolute-value
term is nonnegative.
\end{proof}

\begin{remark}\label{rem:comparison}
Inequality \eqref{eq:asymmetric} is Bre\v{s}ar's bound
\eqref{eq:bresar} plus $\rho(H)\gamma(G)/3$. Inequality
\eqref{eq:symmetric} gives 
$\tfrac23\gamma(G)\gamma(H)$ only when
$\rho(G)/\gamma(G)=\rho(H)/\gamma(H)$; whenever the two factors have different
packing-to-domination ratios, one gains a term proportional to the difference.

On the other hand, \eqref{eq:symmetric} does not contain the classical packing
bound
\[
 \gamma(G\mathbin{\square}H)
 \geq\max\{\rho(G)\gamma(H),\rho(H)\gamma(G)\}.
\]
For example, if $\rho(H)=\gamma(H)$, then \eqref{eq:asymmetric} reads
\[
 \gamma(G\mathbin{\square}H)
 \geq\gamma(G)\gamma(H)-\frac{\rho(G)\gamma(H)}3,
\]
whereas the packing bound gives
$\gamma(G\mathbin{\square}H)\geq\gamma(G)\gamma(H)$. The two bounds should
therefore be used together.
\end{remark}

\section{Limits of the method}\label{sec:limits}

Steiner \cite{Steiner2026} normalizes the relevant parameters as
\[
 x_1:=\frac{\rho(G)}{\gamma(G)},
 \qquad
 x_2:=\frac{\rhoTwo(G)}{2\gamma(G)},
 \qquad
 y_1:=\frac{\rho(H)}{\gamma(H)},
 \qquad
 y_2:=\frac{\rhoTwo(H)}{2\gamma(H)}.
\]
Lemma~\ref{lem:packing-upper} and $\rhoTwo\geq2\rho$ give
$0<x_1\leq x_2\leq1$ and $0<y_1\leq y_2\leq1$.
After division by $\gamma(G)\gamma(H)$, the bounds of this paper and of
\cite{Steiner2026} become lower bounds for
$\gamma(G\mathbin{\square}H)/(\gamma(G)\gamma(H))$ of the following forms,
each together with its mirror image obtained by exchanging
$(x_1,x_2)$ and $(y_1,y_2)$:
\begin{align}
 &\frac{2+y_1-x_1}{3},
 &&\text{by Corollary~\ref{cor:asymmetric}},
 \label{eq:norm18}\\
 &\frac{3+y_1-2x_2+x_2y_1}{4},
 &&\text{by \eqref{eq:packing-function-global}},
 \label{eq:norm19}\\
 &x_1+y_2(1-x_1),
 &&\text{by \cite[Corollary 3.3]{Steiner2026}, from \cite{HouLu2009}},
 \label{eq:norm20}\\
 &\frac{2-x_1}{3},\quad\frac34-\frac{x_2}{2},\quad y_1,\quad\frac12,
 &&\text{by \eqref{eq:bresar}, \eqref{eq:steiner-packing}, the packing bound,
 and \eqref{eq:clark-suen}}.
 \label{eq:norm21}
\end{align}
Steiner's constant \eqref{eq:steiner-main} is the minimum over the box
$0<x_1\leq x_2\leq1$, $0<y_1\leq y_2\leq1$, of the maximum of the six
functions in \eqref{eq:norm20}--\eqref{eq:norm21}, without
\eqref{eq:norm18}--\eqref{eq:norm19}. Adding the new inequalities raises this
minimum to $2/3$ and no further.

\begin{proposition}\label{prop:limits}
Let $\Phi(x_1,x_2,y_1,y_2)$ be the maximum of all functions in
\eqref{eq:norm18}--\eqref{eq:norm21} and their mirror images. Then
\[
 \min\Phi=\frac23
\]
over the box $0<x_1\leq x_2\leq1$, $0<y_1\leq y_2\leq1$. The minimum is
attained at every point
\[
 (x_1,x_2,y_1,y_2)=\left(p,\frac12,p,\frac12\right),
 \qquad 0<p\leq\frac13.
\]
Moreover, at each of these points the linear program obtained from
Lemmas~\ref{lem:sum-L}--\ref{lem:leftover}, with
$\gamma(G)=\gamma(H)=1$, $r=y_1$, $s=1-y_1$, $\rho(G)=x_1$, and
$\rhoTwo(G)=2x_2$, also has optimal value $2/3$.
\end{proposition}

\begin{proof}
By \eqref{eq:norm18} and its mirror image,
\[
 \Phi\geq\frac{2+|x_1-y_1|}{3}\geq\frac23
\]
throughout the box. At $(p,\tfrac12,p,\tfrac12)$, the functions in
\eqref{eq:norm18}--\eqref{eq:norm21} take the values
\[
 \frac23;\qquad
 \frac{2+\tfrac32p}{4}\leq\frac23;\qquad
 \frac{1+p}{2}\leq\frac23;\qquad
 \frac{2-p}{3},\ \frac12,\ p,\ \frac12,
\]
all at most $2/3$ for $0<p\leq1/3$. The mirror images have the same values,
so $\Phi=2/3$ there.

For the linear program, Lemma~\ref{lem:packing-blocks},
\eqref{eq:second-rho}, \eqref{eq:second-rho2}, and
\eqref{eq:leftover-rho}--\eqref{eq:leftover-rho2} combined with
Lemma~\ref{lem:sum-L} give the constraints
\begin{align*}
 a+b&\geq r,\\
 a+\frac b2&\geq r\left(1-\frac{x_1}{2}\right),\\
 a+\frac b3&\geq r\left(1-\frac{2x_2}{3}\right),\\
 \frac32w+\frac b2&\geq s\left(1-\frac{x_1}{2}\right),\\
 \frac43w+\frac b3&\geq s\left(1-\frac{2x_2}{3}\right),
\end{align*}
with $a,b,w\geq0$ and objective $a+b+w$. At
$r=x_1=p$, $s=1-p$, and $x_2=1/2$, put
\[
 a=p-p^2,
 \qquad b=p^2,
 \qquad
 w=\frac23(1-p)\left(1-\frac p2\right)-\frac{p^2}{3}.
\]
Then $a+b=p$ and $a+b/2=p-p^2/2$, so the first two constraints hold with
equality. Also
$a+b/3=p-\tfrac23p^2\geq\tfrac23p$ for $p\leq1/2$. The fourth constraint
holds with equality by the choice of $w$, and the fifth reduces to
$3p^2-6p+2\geq0$, which holds for
$p\leq1-1/\sqrt3$. All variables are nonnegative for $0<p\leq1/3$, and
\[
 a+b+w
 =p+\frac23(1-p)\left(1-\frac p2\right)-\frac{p^2}{3}
 =\frac23.
\]
The optimal value is at least $2/3$ by the proof of
Theorem~\ref{thm:main-packing}, so it equals $2/3$.
\end{proof}

\begin{remark}[The aggregate linear program]\label{rem:lp-exact}
The proof of Theorem~\ref{thm:main-packing} uses only
Lemmas~\ref{lem:sum-L}, \ref{lem:packing-blocks},
\eqref{eq:second-rho}, and \eqref{eq:leftover-rho}. If one minimizes
$a+b+w$ over the corresponding aggregate constraints, including
$a,b,w\geq0$, the optimal value is
\begin{equation}
 \max\left\{
 r\gamma(G),
 \frac{r\gamma(G)}3
 +\frac{(r+s)(2\gamma(G)-\rho(G))}{3}
 \right\}.
 \label{eq:exact-lp}
\end{equation}
Indeed, the second term in \eqref{eq:exact-lp} is the dual lower bound used in
the proof of Theorem~\ref{thm:main-packing}, while the first is
Lemma~\ref{lem:packing-blocks}. When
\begin{equation}
 s(2\gamma(G)-\rho(G))\geq r\rho(G),
 \label{eq:lp-regime}
\end{equation}
the second term is the larger one and an optimal solution is
\[
 a=r(\gamma(G)-\rho(G)),
 \qquad b=r\rho(G),
 \qquad
 w=\frac23s\left(\gamma(G)-\frac{\rho(G)}2\right)-\frac13r\rho(G).
\]
All four aggregate inequalities are then tight. If
\eqref{eq:lp-regime} fails, the exact LP value is $r\gamma(G)$; for example,
one may take $w=0$ and choose
$b=s(2\gamma(G)-\rho(G))\leq r\rho(G)$ and
$a=r\gamma(G)-b$.

At the extremal points used in Proposition~\ref{prop:limits}, condition
\eqref{eq:lp-regime} holds. There the packing parts are negligible,
\eqref{eq:asymmetric} is $\tfrac23\gamma(G)\gamma(H)$ up to lower-order terms,
and the loss comes from applying Lemma~\ref{lem:bresar-steiner} to the sets
$L_i$ of leftover parts, where it gives
$\gamma_G(L_i)\lesssim |L_i|/2$. Lemma~\ref{lem:steiner-subset} improves the
coefficient $1/2$ to $1/3$ but pays $\rhoTwo(G)/3$ per part, which explains why
\eqref{eq:norm19} does not help when $\rhoTwo(G)\approx\gamma(G)$. The
Hou--Lu bound \eqref{eq:norm20} likewise gives only about $1/2$ in this regime.

Thus a further improvement requires an estimate strong when
$\rho(G)\ll\gamma(G)$ and $\rhoTwo(G)\approx\gamma(G)$. In this regime the
half-integral packing furnished by a maximum $2$-packing function has value
$\rhoTwo(G)/2\approx\gamma(G)/2$; this is a lower bound for the fractional
domination number, but it does not by itself determine that parameter. The
hypothetical extremal dominating set described by the linear program covers
only about $\gamma(G)$ of the $|\V(G)|$ fibers in each $G$-layer above a vertex
of $\N_H(p)$; the domination of those layers, expressed by
Lemma~\ref{lem:layer-dominates}, is information not used by the aggregate
argument.
\end{remark}

\begin{remark}\label{rem:free-parameters}
Proposition~\ref{prop:limits} concerns only the optimization problem in which
$x_1,x_2,y_1,y_2$ are treated as free parameters subject to
$x_1\leq x_2$ and $y_1\leq y_2$. It does not assert that
$(\tfrac14,\tfrac12,\tfrac14,\tfrac12)$, or any point with $\Phi=2/3$, is
realized or asymptotically realized by pairs of graphs. Additional universal
relations among $\rho$, $\rhoTwo$, and $\gamma$ would restrict the box and
could raise the certified constant.
\end{remark}

\textbf{Concluding remarks:} The constant $2/3$ of Theorem~\ref{thm:two-thirds} comes from three elementary
ingredients: the packing parts of Lemma~\ref{lem:packing-L-empty}, the layer
bookkeeping of Lemma~\ref{lem:sum-L}, and the second count of
Lemma~\ref{lem:second-count}, combined without first adding the inequalities
of Lemmas~\ref{lem:packing-blocks} and \ref{lem:leftover}. By
Proposition~\ref{prop:limits}, the resulting system of inequalities, even when
strengthened by the bounds of Steiner and of Hou and Lu, certifies
$2/3$ when the normalized graph parameters are treated independently.
Remark~\ref{rem:lp-exact} identifies both the responsible parameter regime and
the information discarded by the aggregate argument.

One natural strengthening of Theorem~\ref{thm:main-packing} is false. The
coefficient of $r$ in \eqref{eq:main-packing} is
$\gamma(G)-\rho(G)/3$, whereas the classical packing bound suggests
$\gamma(G)$. However, the inequality
\[
 \gamma(G\mathbin{\square}H)
 \geq r\gamma(G)+\frac{s(2\gamma(G)-\rho(G))}{3}
\]
fails for $G=K_1$ and $H=P_5$ with $P$ equal to the two end vertices: here
$r=2$, $s=1$, and
$\gamma(G\mathbin{\square}H)=2<2+1/3$. Whether the weaker version with $s$
replaced by $\gamma(H)-\rho(H)$ holds is open; it is implied by Vizing's
conjecture. As the extremal analysis above shows, that inequality alone would
not improve the universal constant.
\bigskip

\noindent \textbf{Declarations}

\medskip
\noindent\textbf{\small Funding:} {\small The research of Mohsen Aliabadi is partially supported by an AMS--Simons Research Enhancement Grant (award no.~SFI-MPS-T-Institutes-00021934).} 

\medskip
\noindent \textbf{\small Conflict of interest:} {\small The authors declare that they have no conflict of interest.}

\medskip
\noindent \textbf{\small Data availability statement:} {\small Data sharing is not applicable to this article as no datasets were generated or analyzed during the current study.}
\medskip

\noindent \textbf{\small AI disclosure:} During the preparation of this work, the authors used Claude (Fable 5) for language editing, to improve the clarity and readability of the manuscript, and for computational assistance related to the proofs of Theorems~\ref{thm:main-packing} and~\ref{thm:packing-function}. The authors subsequently reviewed, verified, and edited all AI-assisted content as needed and take full responsibility for the content of the published article.

\bigskip

\end{document}